\documentclass{amsart}
\usepackage{amsmath,amsthm,amssymb}

\makeatletter
    
    \@addtoreset{equation}{section}
  \makeatother

\newtheorem{definition}{Definition}[section]
\newtheorem{proposition}[definition]{Proposition}
\newtheorem{theorem}[definition]{Theorem}
\newtheorem{lemma}[definition]{Lemma}

\newtheorem{remark}{Remark}[section]

\title[Semilinear damped wave equation]{Critical exponent for the semilinear wave equation with
scale invariant damping}
\author{Yuta WAKASUGI}
\email{y-wakasugi@cr.math.sci.osaka-u.ac.jp}
\address{Department of Mathematics, Graduate School of Science, Osaka University, Toyonaka, Osaka, 560-0043, Japan}
\begin{document}
\begin{abstract}
In this paper we consider the critical exponent problem for the semilinear damped wave equation
with time-dependent coefficients.
We treat the scale invariant cases.
In this case the asymptotic behavior of the solution is very delicate and
the size of coefficient plays an essential role.
We shall prove that if the power of the nonlinearity is greater than the Fujita exponent,
then there exists a unique global solution with small data, provided that
the size of the coefficient is sufficiently large.
We shall also prove some blow-up results
even in the case that the coefficient is sufficiently small.
\end{abstract}

\maketitle
\section{Introduction}
We consider the Cauchy problem for the semilinear damped wave equation
\begin{equation}
\label{eq11}
	\left\{\begin{array}{ll}
	\displaystyle u_{tt}-\Delta u+\frac{\mu}{1+t}u_t=|u|^p,&(t,x)\in (0,\infty)\times\mathbf{R}^n,\\[5pt]
	(u,u_t)(0,x)=(u_0,u_1)(x),&x\in\mathbf{R}^n,
	\end{array}\right.
\end{equation}
where $\mu>0$, $(u_0,u_1)\in H^1(\mathbf{R}^n)\times L^2(\mathbf{R}^n)$ have compact support and
$1<p\le \frac{n}{n-2} (n\ge 3), 1<p<\infty (n=1,2)$.
Our aim is to determine the critical exponent $p_c$, which is the number defined by the following property:

If $p_c<p$, all small data solutions of (\ref{eq11}) are global;
if $1<p\le p_c$, the time-local solution cannot be extended time-globally for some data regardless of smallness.

We note that the linear part of (\ref{eq11}) is invariant with respect to the hyperbolic scaling
$$
	\tilde{u}(t,x)=u(\lambda(1+t)-1,\lambda x).
$$
In this case the asymptotic behavior of solutions is very delicate.
It is known that the size of the damping term $\mu$ plays an essential role.
The damping term $\mu/(1+t)$ is known as the borderline
between the \textit{effective} and \textit{non-effective} dissipation,
here \textit{effective} means that the solution behaves like that of the corresponding
parabolic equation, and
\textit{non-effective} means that
the solution behaves like that of the free wave equation.
Concretely, for linear damped wave equation
\begin{equation}
\label{eq102}
	u_{tt}-\Delta u+(1+t)^{-\beta}u_t=0,
\end{equation}
if $-1<\beta<1$, then the solution $u$ has the same $L^p$-$L^q$ decay rates as those of
the solution of the corresponding heat equation
\begin{equation}
\label{eq103}
	-\Delta v +(1+t)^{-\beta}v_t=0.
\end{equation}
Moreover, if $-1/3<\beta<1$ then the asymptotic profile of $u$ is given by a solution of (\ref{eq103})
(see \cite{W3}).
This is called the \textit{diffusion phenomenon}.
In particular, the constant coefficient case $\beta=0$
has been investigated for a long time.
We refer the reader to \cite{M,N}.
On the other hand, if $\beta>1$ then
the asymptotic profile of the solution of (\ref{eq102}) is given by that of the free wave equation
$\square w=0$ (see \cite{W2}).

Wirth \cite{W1} considered the linear problem
\begin{equation}
\label{eq12}
	\left\{\begin{array}{l}
	\displaystyle u_{tt}-\Delta u+\frac{\mu}{1+t}u_t=0,\\[5pt]
	(u,u_t)(0,x)=(u_0,u_1)(x).
	\end{array}\right.
\end{equation}
He proved several $L^p$-$L^q$ estimates for the solutions to (\ref{eq12}).
For example, if $\mu>1$ it follows that
\begin{align*}
	\|u(t)\|_{L^q}&\lesssim (1+t)^{\max\{ -\frac{n-1}{2}\left(\frac{1}{p}-\frac{1}{q}\right)-\frac{\mu}{2},
		-n\left(\frac{1}{p}-\frac{1}{q}\right)\}}(\|u_0\|_{H^s_p}+\|u_1\|_{H_p^{s-1}}),\\
	\|(u_t, \nabla u)(t)\|_{L^q}&\lesssim (1+t)^{\max\{ -\frac{n-1}{2}\left(\frac{1}{p}-\frac{1}{q}\right)-\frac{\mu}{2},
		-n\left(\frac{1}{p}-\frac{1}{q}\right)-1\}}(\|u_0\|_{H^{s+1}_p}+\|u_1\|_{H_p^{s}}),
\end{align*}
where $1<p\le 2$, $1/p+1/q=1$ and $s=n(1/p-1/q)$.
This shows that if $\mu$ is sufficiently large then
the solution behaves like that of the corresponding heat equation
\begin{equation}
\label{eq13}
	\frac{\mu}{1+t}v_t-\Delta v=0
\end{equation}
as $t\rightarrow \infty$,
and if $\mu$ is sufficiently small then the solution behaves like that of the free wave equation.
We mention that for the wave equation with space-dependent damping
$\square u+V_0\langle x\rangle^{-1}u_t=0$,
a similar asymptotic behavior is obtained by Ikehata, Todorova and Yordanov \cite{ITY}.

The Gauss kernel of (\ref{eq13}) is given by
$$
	G_{\mu}(t,x)=
	\left(\frac{\mu}{2\pi ((1+t)^2-1)}\right)^{\frac{n}{2}}e^{-\frac{\mu|x|^2}{2((1+t)^2-1)}}.
$$
We can obtain the $L^p$-$L^q$ estimates of the solution of (\ref{eq13}).
In fact, it follows that
$$
	\|v(t)\|_{L^q}\lesssim (1+t)^{-n\left(\frac{1}{p}-\frac{1}{q}\right)}\|v(0)\|_{L^p}
$$
for $1\le p\le q\le \infty$.
In particular, taking $q=2$ and $p=1$, we have
$$
	\|v(t)\|_{L^2}\lesssim (1+t)^{-\frac{n}{2}}\|v(0)\|_{L^1}.
$$
From the point of view of the diffusion phenomenon,
we expect that the same type estimate holds for the solution of (\ref{eq12}) when $\mu$ is large.
To state our results, we introduce an auxiliary function
$$
	\psi(t,x):=\frac{a|x|^2}{(1+t)^2},\quad a=\frac{\mu}{2(2+\delta)}
$$
with a positive parameter $\delta$.
We have the following linear estimate:
\begin{theorem}
For any $\varepsilon>0$, there exist constants $\delta>0$ and $\mu_0>1$ such that
for any $\mu\ge\mu_0$ the solution of $(\ref{eq12})$ satisfies
\begin{eqnarray}
\label{eq14}
	\int_{\mathbf{R}^n}e^{2\psi}u^2dx&\le&C_{\mu,\varepsilon}(1+t)^{-n+\varepsilon}\\
\label{eq15}
	\int_{\mathbf{R}^n}e^{2\psi}(u_t^2+|\nabla u|^2)dx
		&\le&C_{\mu,\varepsilon}(1+t)^{-n-2+\varepsilon}
\end{eqnarray}
for $t\ge 0$,
where $C_{\mu,\varepsilon}$ is a positive constant depending on $\mu, \varepsilon$
and $\|(u_0,u_1)\|_{H^1\times L^2}$.
\end{theorem}
\begin{remark}
The constant $\mu_0$ depends on $\varepsilon$.
The relation is
$$
	\mu_0\sim \varepsilon^{-2}.
$$
Therefore, as $\varepsilon$ smaller, $\mu_0$ must be larger.
We can expect that $\varepsilon$ can be removed and the same result holds for
much smaller $\mu$.
However, we have no idea for the proof.
\end{remark}

We also consider the critical exponent problem for (\ref{eq11}).
For the corresponding heat equation (\ref{eq13})
with nonlinear term $|u|^p$,
the critical exponent is given by
$$
	p_F:=1+\frac{2}{n},
$$
which is well known as the Fujita critical exponent (see \cite{F}).
Thus, we can expect that the critical exponent of (\ref{eq11}) is also given by $p_F$
if $\mu$ is sufficiently large.

For the damped wave equation with constant coefficient
$$
	u_{tt}-\Delta u+u_t=|u|^p,
$$
Todorova and Yordanov \cite{TY1}
proved the critical exponent is given by $p_F$.
Lin, Nishihara and Zhai \cite{LNZ} (see also Nishihara\cite{N1})
extended this result to time-dependent coefficient cases
$$
	u_{tt}-\Delta u+(1+t)^{-\beta}u_t=|u|^p
$$
with $-1<\beta<1$.
They proved that $p_F$ is still critical.
Recently, D'abbicco, Lucente and Reissig \cite{DLR} extended
this result to more general effective $b(t)$ by using the linear decay estimates,
which are established by Wirth \cite{W3}.
For the scale-invariant case (\ref{eq11}), very recently, D'abbicco \cite{D} proved the existence of the global solution
with small data for (\ref{eq11}) in the case $n=1,2, \mu\ge n+2$ and $p_F<p$.
He also obtained the decay rates of the solution without any loss $\varepsilon$.

Our main result is following:

\begin{theorem}
Let $p>p_F$ and $0<\varepsilon<2n(p-p_F)/(p-1)$.
Then there exist constants $\delta>0$ and $\mu_0>1$ having the following property:
if $\mu\ge\mu_0$ and
$$
	I_0^2:=\int_{\mathbf{R}^n}e^{2\psi(0,x)}(u_0^2+|\nabla u_0|^2+u_1^2)dx
$$
is sufficiently small, then there exists a unique solution
$u\in C([0,\infty);H^1(\mathbf{R}^n))\cap C^1([0,\infty);L^2(\mathbf{R}^n))$
of $(\ref{eq11})$
satisfying
\begin{eqnarray}
\label{eq16}
	\int_{\mathbf{R}^n}e^{2\psi}u^2dx&\le&C_{\mu,\varepsilon}I_0^2(1+t)^{-n+\varepsilon}\\
\label{eq17}
	\int_{\mathbf{R}^n}e^{2\psi}(u_t^2+|\nabla u|^2)dx
		&\le&C_{\mu,\varepsilon}I_0^2(1+t)^{-n-2+\varepsilon}
\end{eqnarray}
for $t\ge 0$,
where $C_{\mu,\varepsilon}$ is a positive constant depending on $\mu$ and $\varepsilon$.
\end{theorem}

\begin{remark}
Similarly as before, we note that $\mu_0$ depends on $\varepsilon$.
The relation is
$$
	\mu_0\sim \varepsilon^{-2}\sim (p-p_F)^{-2}.
$$
Therefore, as $p$ is closer to $p_F$, $\mu_0$ must be larger.
Thus, we can expect that $\varepsilon$
can be removed and the same result holds for much smaller $\mu$.
As mentioned above,
D'abbicco \cite{D} obtained an affirmative result for this expectation when $n=1,2$.
However, we have no idea for the higher dimensional cases.
\end{remark}

We prove Theorem 1.2 by a multiplier method
which is essentially developed in \cite{TY1}.
Lin, Nishihara and Zhai \cite{LNZ} refined this method
to fit the damping term $b(t)=(1+t)^{-\beta}$ with $-1<\beta<1$.
They used the property $\beta<1$
and so we cannot apply their method directly to our problem (\ref{eq11}).
Therefore, we need a further modification.
Instead of the property $\beta<1$, we assume that $\mu$ is sufficiently large
and modify the parameters used in the calculation.

\begin{remark}
We can also treat other nonlinear terms, for example $-|u|^p, \pm |u|^{p-1}u$.
\end{remark}

We also have a blow-up result when $\mu>1$ and $1<p\le p_F$.

\begin{theorem}
Let $1<p\le p_F$ and $\mu>1$.
Moreover, we assume that
$$
	\int_{\mathbf{R}^n}(\mu-1)u_0+u_1dx>0.
$$
Then there is no global solution for $(\ref{eq11})$.
\end{theorem}
\begin{remark}
Theorem 1.3 is essentially included in a recent work by D'abbicco and Lucente \cite{DL}.
In this paper we shall give a much simpler proof.
\end{remark}

One of our novelty is blow-up results for the non-effective damping cases.
We also obtain blow-up results in the case $0<\mu\le 1$.
\begin{theorem}
Let $0<\mu\le 1$ and
$$
	1<p\le 1+\frac{2}{n+(\mu-1)}.
$$
We also assume
$$
	\int_{\mathbf{R}^n}u_1(x)dx>0.
$$
Then there is no global solution for $(\ref{eq11})$.
\end{theorem}

\begin{remark}
In Theorem 1.4, we do not put any assumption on the data $u_0$,
and the blow-up results hold even for the case $p\ge p_F$.
We can interpret this phenomena as that
the equation $(\ref{eq11})$ loses the parabolic structure
and recover the hyperbolic structure if $\mu$ is sufficiently small.
\end{remark}

We prove this theorem by a test-function method developed by Qi S. Zhang \cite{Z}.
In the same way of the proof of Theorem 1.4,
we can treat the damping terms $(1+t)^{-\beta}$ with $\beta>1$ (see Remark 3.1).

In the next section, we give a proof of Theorem 1.2.
We can prove Theorem 1.1 by the almost same way,
and so we omit the proof.
In Section 3, we shall prove Theorem 1.3 and Theorem 1.4.

\section{Proof of Theorem 1.2}
In this section we prove our main result.
First, we prepare some notation and terminology.
We put
$$
	\|f\|_{L^p(\mathbf{R}^n)}:=
	\left(\int_{\mathbf{R}^n}|f(x)|^pdx\right)^{1/p}.
$$
By $H^1(\mathbf{R}^n)$ we denote the usual Sobolev space.
For an interval $I$ and a Banach space $X$, we define $C^r(I;X)$ as the Banach space whose
element is an $r$-times continuously differentiable mapping from $I$ to $X$ with respect to the topology in $X$.
The letter $C$ indicates the generic constant, which may change from a line to the next line. 
We also use the symbols $\lesssim$ and $\sim$.
The relation $f\lesssim g$ means $f\le Cg$ with some constant $C>0$
and $f\sim g$ means $f\lesssim g$ and $g\lesssim f$.

We first describe the local existence:

\begin{proposition}\label{prop21}
For any $p>1, \mu>0$ and $\delta>0$, there exists $T_m\in (0,+\infty]$ depending on $I_0^2$ such that
the Cauchy problem $(\ref{eq11})$
has a unique solution $u\in C([0,T_m);H^1(\mathbf{R}^n))\cap C^1([0,T_m);L^2(\mathbf{R}^n))$,
and if $T_m <+\infty$ then we have
$$
	\liminf_{t\to T_m}\int_{\mathbf{R}^n}
	e^{\psi(t,x)}
		(u_t^2+|\nabla u|^2+u^2)dx=+\infty.$$
\end{proposition}
We can prove this proposition by standard arguments (see \cite{IT}).
We prove a priori estimate for the following functional:

$$
	M(t):=\sup_{0\le\tau\le t}
	\left\{(1+\tau)^{n+2-\varepsilon}\int_{\mathbf{R}^n}e^{2\psi}(u_t^2+|\nabla u|^2)dx
		+(1+\tau)^{n-\varepsilon}\int_{\mathbf{R}^n}e^{2\psi}u^2dx\right\}.
$$

We put $b(t)=\frac{\mu}{1+t}$ and $f(u)=|u|^{p}$.
By a simple calculation, we have
$$
	-\psi_t=\frac{2}{1+t}\psi,\quad \nabla\psi=\frac{2ax}{(1+t)^2},\quad
	\frac{|\nabla\psi|^2}{-\psi_t}=\frac{b(t)}{2+\delta}
$$
and
$$
	\Delta\psi=\frac{n}{2+\delta}\frac{b(t)}{1+t}=:\left(\frac{n}{2}-\delta_1\right)\frac{b(t)}{1+t}.
$$
Here and after,
$\delta_i(i=1,2,\ldots)$
denote a positive constant depending only on
$\delta$
such that
$$
	\delta_i\rightarrow 0^+ \quad \mbox{as} \quad\delta\rightarrow 0^+.
$$
Multiplying (\ref{eq11}) by $e^{2\psi}u_t$, we obtain
\begin{eqnarray}
	\lefteqn{\frac{\partial}{\partial t}\left[ \frac{e^{2\psi}}{2}(u_t^2+|\nabla u|^2)\right]
		-\nabla\cdot(e^{2\psi}u_t\nabla u)}\label{eq21}\\
	&&+e^{2\psi}\left(b(t)-\frac{|\nabla\psi|^2}{-\psi_t}-\psi_t\right)u_t^2
		+\underbrace{\frac{e^{2\psi}}{-\psi_t}|\psi_t\nabla u-u_t\nabla\psi|^2}_{T_1}\nonumber\\
	&&=\frac{\partial}{\partial t}\left[e^{2\psi}F(u)\right]+2e^{2\psi}(-\psi_t)F(u),\nonumber
\end{eqnarray}
where $F$ is the primitive of $f$ satisfying $F(0)=0$.
Using the Schwarz inequality, we can calculate
$$
	T_1\ge e^{2\psi}
		\left(\frac{1}{5}(-\psi_t)|\nabla u|^2-\frac{b(t)}{4(2+\delta)}u_t^2\right).
$$
From this and integrating (\ref{eq21}), we have
\begin{align}
\label{eq22}
	\frac{d}{dt}\int_{\mathbf{R}^n}\frac{e^{2\psi}}{2}(u_t^2+|\nabla u|^2)dx
		&+\int_{\mathbf{R}^n}e^{2\psi}\left\{\left(\frac{b(t)}{4}-\psi_t\right)u_t^2+\frac{-\psi_t}{5}|\nabla u|^2
		\right\}dx\\
	&\le\frac{d}{dt}\int_{\mathbf{R}^n}e^{2\psi}F(u)dx+2e^{2\psi}(-\psi_t)F(u)dx.\nonumber
\end{align}
On the other hand, by multiplying (\ref{eq11}) by $e^{2\psi}u$, it follows that
\begin{eqnarray*}
\lefteqn{\frac{\partial}{\partial t}\left[e^{2\psi}\left(uu_t+\frac{b(t)}{2}u^2\right)\right]
		-\nabla\cdot (e^{2\psi}u\nabla u)}\\
&&+e^{2\psi}\Big\{|\nabla u|^2+\left(-\psi_t+\frac{1}{2(1+t)}\right)b(t)u^2
	+\underbrace{2u\nabla\psi\cdot\nabla u}_{T_2}-2\psi_tuu_t-u_t^2\Big\}\nonumber\\
&&=e^{2\psi}uf(u).\nonumber
\end{eqnarray*}
We calculate
\begin{eqnarray*}
	e^{2\psi}T_2&=&4e^{2\psi}u\nabla\psi\cdot\nabla u-2e^{2\psi}u\nabla\psi\cdot\nabla u\\
	&=&4e^{2\psi}u\nabla\psi\cdot\nabla u
		-\nabla\cdot(e^{2\psi}u^2\nabla\psi)+2e^{2\psi}u^2|\nabla\psi|^2+e^{2\psi}(\Delta\psi)u^2
\end{eqnarray*}
and have
\begin{eqnarray}
\label{eq23}
	\lefteqn{\frac{\partial}{\partial t}\left[e^{2\psi}\left(uu_t+\frac{b(t)}{2}u^2\right)\right]
		-\nabla\cdot(e^{2\psi}(u\nabla u+u^2\nabla\psi))}\\
	&&+e^{2\psi}\Big\{\underbrace{|\nabla u|^2+4u\nabla u\cdot\nabla\psi
		+((-\psi_t)b(t)+2|\nabla\psi|^2)u^2}_{T_3}\nonumber\\
	&&+(n+1-2\delta_1)\frac{b(t)}{2(1+t)}u^2-2\psi_tuu_t-u_t^2\Big\} =e^{2\psi}uf(u).\nonumber
\end{eqnarray}
$T_3$ is estimated as
\begin{eqnarray*}
	T_3&=&\left(1-\frac{4}{4+\delta/2}\right)|\nabla u|^2+\frac{\delta}{2}|\nabla\psi|^2u^2
		+\left|\frac{2}{\sqrt{4+\delta/2}}\nabla u+\sqrt{4+\delta/2}\nabla\psi\right|^2\\
	&\ge&\delta_2(|\nabla u|^2+b(t)(-\psi_t)u^2).
\end{eqnarray*}
Thus, we can rewrite (\ref{eq23}) as
\begin{eqnarray*}
	\lefteqn{\frac{\partial}{\partial t}\left[e^{2\psi}\left(uu_t+\frac{b(t)}{2}u^2\right)\right]
		-\nabla\cdot (e^{2\psi}(u\nabla u+u^2\nabla\psi))}\\
	&&+e^{2\psi}\left\{
		\delta_2(|\nabla u|^2+b(t)(-\psi_t)u^2)+(n+1-2\delta_2)\frac{b(t)}{2(1+t)}u^2
			-2\psi_tuu_t-u_t^2\right\}\nonumber\\
	&&\le e^{2\psi}uf(u).\nonumber
\end{eqnarray*}
Integrating the above inequality and then
multiplying by a large parameter $\nu$ and adding $(1+t)\times$ (\ref{eq22}), we obtain
\begin{eqnarray*}
	\lefteqn{
		\frac{d}{dt}\left[
			\int e^{2\psi}\left\{
				\frac{1+t}{2}(u_t^2+|\nabla u|^2)+\nu uu_t+\frac{\nu b(t)}{2}u^2
			\right\}dx
		\right]
	}\\
	&&+\int e^{2\psi}
		\Big\{
			\Big(
				\underbrace{\frac{\mu}{4}-\nu-\frac{1}{2}}_{T_4}+(-\psi_t)(1+t)
			\Big)u_t^2
			+\Big(
				\underbrace{\nu\delta_2-\frac{1}{2}}_{T_5}+\frac{(-\psi_t)(1+t)}{5}
			\Big)|\nabla u|^2\\
	&&+\nu\delta_2b(t)(-\psi_t)u^2+(n+1-2\delta_1)\frac{\nu b(t)}{2(1+t)}u^2
		+\underbrace{2\nu(-\psi_t)uu_t}_{T_6}\Big\}dx\\
	&&\le \frac{d}{dt}
		\left[
			(1+t)\int e^{2\psi}F(u)dx
		\right]
		+C\int e^{2\psi}(1+(1+t)(-\psi_t))|u|^{p+1}dx.
\end{eqnarray*}
We put the condition for $\mu$ and $\nu$ as
\begin{eqnarray}
\label{eq24}
	\frac{\mu}{4}-\nu-\frac{1}{2}&>&0\\
\label{eq25}
	\nu\delta_2-\frac{1}{2}&>&0.
\end{eqnarray}
Then the terms $T_4$ and $T_5$ are positive.
Using the Schwarz inequality, we can estimate $T_6$ as
$$
	|T_6|\le
		\frac{1}{2}(-\psi_t)(1+t)u_t^2+\frac{2\nu^2}{1+t}(-\psi_t)u^2.
$$
Now we put an another condition
\begin{equation}
\label{eq26}
	\mu\ge \frac{2\nu}{\delta_2}.
\end{equation}
Then we obtain the following estimate.
\begin{eqnarray}
\label{eq27}
	\lefteqn{
	\frac{d}{dt}\hat{E}(t)+H(t)
		+(n+1-2\delta_1)\frac{\nu b(t)}{2(1+t)}J(t;u^2)}\\
	&\le&
	\frac{d}{dt}[(1+t)J(t;F(u))]+C(J(t;|u|^{p+1})+(1+t)J_{\psi}(t;|u|^{p+1})),\nonumber
\end{eqnarray}
where
$$
	\hat{E}(t):=\int e^{2\psi}\left\{
				\frac{1+t}{2}(u_t^2+|\nabla u|^2)+\nu uu_t+\frac{\nu b(t)}{2}u^2
			\right\}dx,
$$
\begin{align*}
	&H(t)=\int e^{2\psi}\left\{\left(\frac{\mu}{4}-\nu-\frac{1}{2}\right)u_t^2
		+\left(\nu\delta_2-\frac{1}{2}\right)|\nabla u|^2\right\}dx&\\
	&J(t;u)=\int e^{2\psi}udx,\quad J_{\psi}(t;u)=\int e^{2\psi}(-\psi_t)udx.&
\end{align*}
Multiplying (\ref{eq27}) by $(1+t)^{n+1-\varepsilon}$, we have
\begin{eqnarray*}
	\lefteqn{
	\frac{d}{dt}[(1+t)^{n+1-\varepsilon}\hat{E}(t)]
		-\underbrace{(n+1-\varepsilon)(1+t)^{n-\varepsilon}\hat{E}(t)}_{T_7}}\\
	&&+(1+t)^{n+1-\varepsilon}H(t)
		+(n+1-2\delta_1)(1+t)^{n+1-\varepsilon}\frac{\nu b(t)}{2(1+t)}J(t;u^2)\\
	&\le&
	\frac{d}{dt}[(1+t)^{n+2-\varepsilon}J(t;F(u))]\\
		&&+C((1+t)^{n+1-\varepsilon}(J(t;|u|^{p+1})+(1+t)J_{\psi}(t;|u|^{p+1})).\nonumber
\end{eqnarray*}
Now we estimate the bad term $T_7$.
First, by the Schwarz inequality, one can obtain
$$
	|\nu uu_t|\le\frac{\nu}{4\delta_3 b(t)}u_t^2+\delta_3\nu b(t)u^2,
$$
where $\delta_3$ determined later.
From this, $T_7$ is estimated as
\begin{eqnarray*}
	T_7&\le&(n+1-\varepsilon)(1+t)^{n-\varepsilon}\\
		&&\times\int e^{2\psi}\left\{
			\left(\frac{1+t}{2}+\frac{\nu(1+t)}{4\delta_3\mu}\right)u_t^2
			+\frac{1+t}{2}|\nabla u|^2+\frac{\nu b(t)}{2}(1+2\delta_3)u^2
		\right\}dx.
\end{eqnarray*}
We strengthen the conditions (\ref{eq24}) and (\ref{eq25}) as
\begin{eqnarray}
\label{eq28}
	\frac{\mu}{4}-\nu-\frac{1}{2}-(n+1-\varepsilon)\left(\frac{1}{2}+\frac{\nu}{4\delta_3\mu}\right)&>&0,\\
\label{eq29}
	\nu\delta_2-\frac{1}{2}(n+2-\varepsilon)&>&0.
\end{eqnarray}
Moreover, we take $\varepsilon=3\delta_1$ and then choose $\delta_3$ such that
$$
	(n+1-2\delta_1)-(n+1-3\delta_1)(1+2\delta_3)>0.
$$
Under these conditions, we can estimate $T_7$ and obtain
\begin{eqnarray*}
	\frac{d}{dt}[(1+t)^{n+1-\varepsilon}\hat{E}(t)]
	&\le&
	\frac{d}{dt}[(1+t)^{n+2-\varepsilon}J(t;F(u))]\\
	&&+C(1+t)^{n+1-\varepsilon}(J(t;|u|^{p+1})+(1+t)J_{\psi}(t;|u|^{p+1})).
\end{eqnarray*}
By integrating the above inequality, it follows that
\begin{align*}
	(1+t)^{n+1-\varepsilon}\hat{E}(t)
	&\le
	CI_0^2+(1+t)^{n+2-\varepsilon}J(t;|u|^{p+1})\\
	&+C\int_0^t(1+\tau)^{n+1-\varepsilon}(J(\tau;|u|^{p+1})+(1+t)J_{\psi}(\tau;|u|^{p+1}))d\tau.
\end{align*}
By a simple calculation, we have
$$
	(1+t)E(t)+\frac{1}{1+t}J(t;u^2)\le C\hat{E}(t),
$$
where
$$
	E(t):=\int e^{2\psi}(u_t^2+|\nabla u|^2)dx.
$$
Thus, we obtain
\begin{eqnarray}
\label{eq210}
	\lefteqn{(1+t)^{n+2-\varepsilon}E(t)+(1+t)^{n-\varepsilon}J(t;u^2)}\\
	&\le&
	CI_0^2+(1+t)^{n+2-\varepsilon}J(t;|u|^{p+1})\nonumber\\
	&&+C\int_0^t(1+\tau)^{n+1-\varepsilon}(J(\tau;|u|^{p+1})+(1+t)J_{\psi}(\tau;|u|^{p+1}))d\tau.\nonumber
\end{eqnarray}
Now we turn to estimate the nonlinear terms.
We need the following lemma:
\begin{lemma}[Gagliardo-Nirenberg]\label{lem22}
Let $p,q,r \in [1,\infty]$ and $\sigma\in [0,1]$ satisfy
$$\frac{1}{p}=\sigma\left(\frac{1}{r}-\frac{1}{n}\right)+(1-\sigma)\frac{1}{q}$$
except for $p=\infty$ or $r=n$ when $n\ge 2$. Then for some constant $C=C(p,q,r,n)>0$,
the inequality
$$\|h\|_{L^p}\le C\|h\|_{L^q}^{1-\sigma}\|\nabla h\|_{L^r}^{\sigma}
		\quad \mbox{for any}\ h\in C^1_0(\mathbf{R}^n)$$
holds.
\end{lemma}
Noting that
$$
	J(t;|u|^{p+1})=\int \left| e^{\frac{2}{p+1}\psi}u\right|^{p+1}dx
$$
and $\nabla (e^{\frac{2}{p+1}\psi}u)=\frac{2}{p+1}e^{\frac{2}{p+1}\psi}(\nabla\psi)u+
e^{\frac{2}{p+1}\psi}\nabla u$, we apply the above lemma to $J(t;|u|^{p+1})$ with
$\sigma =\frac{n(p-1)}{2(p+1)}$ and have
\begin{align*}
	J(t;|u|^{p+1})&\le
	C\left(\int e^{\frac{4}{p+1}\psi}u^2dx\right)^{\frac{1-\sigma}{2}(p+1)}\\
	&\times \left(\int e^{\frac{4}{p+1}\psi}|\nabla\psi|^2u^2dx
		+\int e^{\frac{4}{p+1}\psi}|\nabla u|^2dx\right)^{\frac{\sigma}{2}(p+1)}.
\end{align*}
We note that
$$
	e^{\frac{4}{p+1}\psi}|\nabla\psi|^2=
	\frac{4a^2|x|^2}{(1+t)^4}e^{\frac{4}{p+1}\psi}\\
	\le C\frac{1}{(1+t)^2}e^{2\psi}
$$
and obtain
\begin{align*}
	J(t;|u|^{p+1})&\le
	C\left(\int e^{2\psi}u^2dx\right)^{\frac{1-\sigma}{2}(p+1)}\\
	&\times \left(\frac{1}{(1+t)^2}\int e^{2\psi}u^2dx
		+\int e^{2\psi}|\nabla u|^2dx\right)^{\frac{\sigma}{2}(p+1)}.
\end{align*}
Therefore, we can estimate
\begin{align*}
	(1+t)^{n+2-\varepsilon}J(t;|u|^{p+1})
	&\le
	(1+t)^{n+2-\varepsilon}
	\{(1+t)^{-(n-\varepsilon)}M(t)\}^{\frac{1-\sigma}{2}(p+1)}\\
	&\times\{(1+t)^{-(n+2-\varepsilon)}M(t)\}^{\frac{\sigma}{2}(p+1)}.
\end{align*}
By a simple calculation, if
\begin{equation}
\label{eq211}
	\varepsilon<\frac{2n\left( p-\left(1+\frac{2}{n}\right)\right)}{p-1}
\end{equation}
then we have
$$
	(1+t)^{n+2-\varepsilon}J(t;|u|^{p+1})
	\le CM(t)^{p+1}.
$$
We note that the conditions (\ref{eq26}), (\ref{eq28}), (\ref{eq29}), (\ref{eq211}) are fulfilled
by the determination of the parameters in the order
$$
	p\rightarrow \varepsilon \rightarrow \delta \rightarrow \nu \rightarrow \mu.
$$
In a similar way, we can estimate the other nonlinear terms.
Consequently, we obtain the a priori estimate
$$
	M(t)\le CI_0^2+CM(t)^{p+1}.
$$
This proves Theorem 1.2.

\section{Proof of Theorem 1.3 and Theorem 1.4}
In this section we first give a proof of Theorem 1.3.
We use a method by Lin, Nishihara and Zhai \cite{LNZ} to transform (\ref{eq11})
into divergence form and then a test-function method by Qi S. Zhang \cite{Z}.

Let $\mu>1$. We
multiply (\ref{eq11}) by a positive function $g(t)\in
C^2([0,\infty))$ and obtain
$$
    (gu)_{tt}-\Delta(gu)-(g^{\prime}u)_t+(-g^{\prime}+gb)u_t=g|u|^p.
$$
We now choose $g(t)$ as the solution of the initial value problem
for the ordinary differential equation
\begin{equation}
\label{eq30}
    \left\{\begin{array}{l}
        -g^{\prime}(t)+g(t)b(t)=1,\quad t>0,\\[5pt]
        g(0)=\frac{1}{\mu-1}.
    \end{array}\right.
\end{equation}
The solution $g(t)$ is explicitly given by
$$
    g(t)=\frac{1}{\mu-1}(1+t).
$$
Thus, we obtain the equation in divergence form
\begin{equation}
\label{eq31}
    (gu)_{tt}-\Delta(gu)-(g^{\prime}u)_t+u_t=g|u|^p.
\end{equation}

Next, we apply a test function method.
We first introduce test functions having the
following properties:
\begin{align}
\tag{t1}
    \eta(t)\in C_0^{\infty}([0,\infty)),\quad 0\le\eta(t)\le 1,\quad
    \eta(t)=\left\{\begin{array}{ll}
        1,&0\le t\le\frac{1}{2},\\
        0,&t\ge 1,
    \end{array}\right.\\
\tag{t2}
    \phi(x)\in C_0^{\infty}(\mathbf{R}^n),\quad0\le\phi(x)\le 1,\quad
    \phi(x)=\left\{\begin{array}{ll}
        1,&|x|\le\frac{1}{2},\\
        0,&|x|\ge 1,
    \end{array}\right.\\
\tag{t3}
    \frac{\eta^{\prime}(t)^2}{\eta(t)}\le C\quad\Big(\frac{1}{2}\le t\le 1\Big),\quad
    \frac{|\nabla\phi(x)|^2}{\phi(x)}\le C\quad \Big(\frac{1}{2}\le |x|\le 1\Big).&\
\end{align}
Let $R$ be a large parameter in $(0,\infty)$. We define the test
function
$$
    \psi_R(t,x):=\eta_R(t)\phi_R(x):=\eta\left(\frac{t}{R}\right)\phi\left(\frac{x}{R}\right),
$$
Let $q$ be the dual of $p$, that is $q=\frac{p}{p-1}$.
Suppose that $u$ is a global solution with initial data $(u_0,u_1)$ satisfying
$$
	\int_{\mathbf{R}^n}((\mu-1)u_0+u_1)dx>0.
$$
We define
$$
    I_R:=\int_{Q_R}g(t)|u(t,x)|^p\psi_R(t,x)^qdxdt,
$$
where $Q_R:=[0,R]\times B_R$ and $B_R:=\{x\in\mathbf{R}^n ; |x|\le R\}$.
By the equation (\ref{eq31}) and integration by parts
one can calculate
\begin{eqnarray*}
    \lefteqn{I_R=-\int_{B_R}((\mu-1)u_0+u_1)\phi^q_Rdx}\\
        &&+\int_{Q_R}gu\partial_t^2(\psi_R^q)dxdt
        +\int_{Q_R}(g^{\prime}u-u)\partial_t(\psi_R^q)dxdt
        -\int_{Q_R}gu\Delta(\psi_R^q)dxdt\\
        &=:&-\int_{B_R}((\mu-1)u_0+u_1)\phi^q_Rdx+J_1+J_2+J_3.
\end{eqnarray*}
By the assumption on the data $(u_0,u_1)$ it follows that
$$
    I_R<J_1+J_2+J_3
$$
for large $R$.
We shall estimate $J_1, J_2$ and $J_3$, respectively.
We use the notation
$$
    \hat{Q}_{R}:=[R/2, R]\times B_R,\quad \tilde{Q}_{R}:=[0,R]\times (B_R\setminus B_{R/2}).
$$
We first estimate $J_3$. By the conditions {\rm (t1)-(t3)} and the H\"older inequality we have the
following estimate:
\begin{eqnarray*}
    |J_3|&\lesssim&R^{-2}\int_{\tilde{Q}_{R}}g(t)|u|\psi_R^{q-1}dxdt \\
    &\lesssim&R^{-2}\left(\int_{\tilde{Q}_{R}}g(t)|u|^p\psi_R^q(t,x)dxdt\right)^{1/p}
        \left(\int_{\tilde{Q}_{R}}g(t)dxdt\right)^{1/q}\\
    &\lesssim&\tilde{I}_{R}^{1/p}R^{\frac{n+2}{q}-2},
\end{eqnarray*}
where
$$
    \tilde{I}_{R}:=\int_{\tilde{Q}_{R}}g(t)|u|^p\psi_R^q(t,x)dxdt.
$$
In a similar way, we can estimate $J_1$ and $J_2$ as
$$
	|J_1|+|J_2|\lesssim \hat{I}_{R}^{1/p}R^{\frac{n+2}{q}-2},\quad
	\hat{I}_{R}:=\int_{\hat{Q}_{R}}g(t)|u|^p\psi_R^q(t,x)dxdt.
$$
Hence, we obtain
\begin{equation}
\label{eq32}
	I_R\lesssim (\tilde{I}_{R}^{1/p}+\hat{I}_{R}^{1/p})R^{\frac{n+2}{q}-2},
\end{equation}
in particular $I_R^{1-1/p}\lesssim R^{\frac{n+2}{q}-2}$.
If $1<p<p_F$, by letting $R\rightarrow \infty$ we have $I_R\rightarrow 0$ and hence
$u=0$, which contradicts the assumption on the data.
If $p=p_F$, we have only
$I_R\le C$ with some constant $C$ independent of $R$.
This implies that $g(t)|u|^p$ is integrable on $(0,\infty)\times\mathbf{R}^n$ and hence
$$
	\lim_{R\to\infty}(\tilde{I}_{R}+\hat{I}_{R})=0.
$$
By (\ref{eq32}), we obtain $\lim_{R\to\infty}I_R=0$.
Therefore, $u$ must be $0$. This also leads a contradiction.

\begin{proof}[Proof of Theorem 1.4]
The proof is almost same as above.
Let $0<\mu\le 1$.
Instead of (\ref{eq30}),
we consider the ordinary differential equation
\begin{equation}
\label{eq34}
        -g^{\prime}(t)+g(t)b(t)=0
\end{equation}
with $g(0)>0$.
We can easily solve this and have
$$
	g(t)=g(0)(1+t)^{\mu}.
$$
Then we have
\begin{equation}
\label{eq35}
	(gu)_{tt}-\Delta (gu)-(g^{\prime}u)_t=g|u|^p.
\end{equation}
Using the same test function $\psi_R(t,x)$ as above, we can calculate
\begin{align*}
	I_R&:=\int_{Q_R}g(t)|u|^p\psi_R^qdxdt\\
	&=-\int_{B_R}g(0)u_1\phi_R^qdx+\sum_{k=1}^3J_k,
\end{align*}
where
$$
	J_1=\int_{Q_R}gu\partial_t^2(\psi_R^q)dxdt,\quad
	J_2=\int_{Q_R}g^{\prime}u\partial_t(\psi_R^q)dxdt,\quad
	J_3=-\int_{Q_R}gu\Delta(\psi_R^q)dxdt.
$$
We note that the term of $u_0$ vanishes and so we put
the assumption only for $u_1$.
We first estimate $J_2$.
Noting $g^{\prime}(t)=\mu g(0)(1+t)^{\mu-1}$, we have
$$
	|J_2|\lesssim \frac{1}{R}\int_{\hat{Q}_R}(1+t)^{\mu-1}|u|\psi_R^{q-1}dxdt.
$$
By noting that
$(1+t)^{\mu-1}\sim g(t)^{1/p}(1+t)^{\mu/q-1}$
and the H\"older inequality, it follows that
\begin{align*}
	|J_2|&\lesssim \frac{1}{R}\left(\int_{\hat{Q}_R}g|u|^p\psi_R^qdxdt\right)^{1/p}
		\left(\int_{R/2}^R\int_{B_R}(1+t)^{\mu-q}dxdt\right)^{1/q}\\
	&\lesssim \hat{I}_R^{1/p}R^{-1+(n+(\mu-q+1))/q},
\end{align*}
where $\hat{I}_R$ is defined as before.
A simple calculation shows
$$
	-1+(n+(\mu-q+1))/q\le 0\ \Leftrightarrow \ p\le 1+\frac{2}{n+(\mu-1)}.
$$
In the same way, we can estimate $J_1$ and $J_3$ as
$$
	|J_1|+|J_3|\lesssim (\hat{I}_R^{1/p}+\tilde{I}_{R}^{1/p})R^{-2+(n+\mu+1)/q},
$$
where $\tilde{I}_R$ is same as before.
It is also easy to see that
$$
	-2+(n+\mu+1)/q\le 0\ \Leftrightarrow \ p\le 1+\frac{2}{n+(\mu-1)}.
$$
Finally, we have
$$
	I_R\lesssim \tilde{I}_R^{1/p}+\hat{I}_R^{1/p}
$$
if $p\le 1+2/(n+(\mu-1))$.
The rest of the proof is same as before.
\end{proof}
\begin{remark}
We can apply the proof of Theorem 1.4 to the wave equation with non-effective damping terms
$$
	\left\{\begin{array}{l}
	u_{tt}-\Delta u+b(t)u_t=|u|^p,\\
	(u,u_t)(0,x)=(u_0,u_1)(x),
	\end{array}\right.
$$
where
$$
	b(t)=(1+t)^{-\beta}
$$
with $\beta>1$.
We can easily solve $(\ref{eq34})$ and have
$$
	g(t)=g(0)\exp\left(\frac{1}{-\beta+1}((1+t)^{-\beta+1}-1)\right).
$$
We note that $g(t)\sim 1$.
The same argument implies that if
$$
	1<p\le 1+\frac{2}{n-1},\quad \int u_1dx>0,
$$
then there is no global solution.
We note that the exponent $1+2/(n-1)$ is greater than the Fujita exponent.
This shows that when $\beta>1$,
the equation loses the parabolic structure even in the nonlinear cases.
One can expect that the critical exponent is given by
the well-known Strauss critical exponent.
However, this problem is completely open as far as the author's knowledge.
\end{remark}

\section*{Acknowledgement}
The author has generous support from Professors Tatsuo Nishitani and Kenji Nishihara.
In particular, Prof. Nishihara gave the author an essential idea for the proof of Theorem 1.4.

\end{document}